\documentclass[twoside,12pt]{article}
\usepackage{amsmath,amssymb}
\date{}
\newtheorem{Lemma}{LEMMA}[section]
\newtheorem{Corollary}[Lemma]{COROLLARY}
\newtheorem{Theorem}[Lemma]{THEOREM}
\newtheorem{Proposition}[Lemma]{PROPOSITION}

\newcommand{\bnum}{\begin{enumerate}}
\newcommand{\enum}{\end{enumerate}}
\newcommand{\bi}{\begin{itemize}}
\newcommand{\ei}{\end{itemize}}
\newcommand{\btab}{\begin{tabular}}
\newcommand{\etab}{\end{tabular}}
\newcommand{\beq}{\begin{eqnarray*}}
\newcommand{\eeq}{\end{eqnarray*}}
\newcommand{\beqn}{\begin{eqnarray}}
\newcommand{\eeqn}{\end{eqnarray}}

\newcommand{\bq}{\begin{equation}}
\newcommand{\eq}{\end{equation}}

\newcommand{\CC}{{\mathcal C}}

\newcommand{\CL}{{\mathcal L}}

\newcommand{\CS}{{\mathcal S}}

\def\phi{\varphi}
\def\epsilon{\varepsilon}

\newcommand{\BR}{\mathbb R}
\newcommand{\BC}{\mathbb C}

\newcommand{\BH}{\mathbb H}
\newcommand{\DP}{{{\rm PG}(3,\BR)}}
\newcommand{\SO}{{{\rm SO}_3 \BR}}
\newcommand{\So}{{{\rm SO}_2 \BR}}

\newcommand{\kasten}{\vbox{\hrule height 8pt width 8.6pt depth -7.4pt
    \hbox{\vrule width 0.6pt height 7.4pt
    \kern 7.4pt \vrule width 0.6pt height 7.4pt}
    \hrule height 0.6pt width 8.6pt}}
\newcommand{\ok}{\hfill\kasten}
\newcommand{\bpf}{\begin{Proof}}
\newcommand{\epf}{\ok\end{Proof}\bigskip\noindent}
\newcommand{\bthm}{\begin{Theorem}}
\newcommand{\ethm}{\end{Theorem}}
\newcommand{\ble}{\begin{Lemma}}
\newcommand{\ele}{\end{Lemma}}
\newcommand{\bprop}{\begin{Proposition}}
\newcommand{\eprop}{\end{Proposition}}
\newcommand{\bcor}{\begin{Corollary}}
\newcommand{\ecor}{\end{Corollary}}

\begin{document}
\title{Parallelisms of ${\rm PG}(3,\mathbb R)$ admitting a 3-dimensional group}

\author{Rainer L\"owen}

\maketitle
\thispagestyle{empty}

\begin{abstract}
Betten and Riesinger \cite{nonreg} constructed Parallelisms of $\DP$ with automorphism group $\SO$ by 
applying the reducible $\SO$-action to a  rotational Betten spread. 
This was generalized in \cite{reduzibel} so as to include oriented parallelisms (i.e., parallelisms of oriented lines). 
In this way, a much larger class of examples was produced. Here we show that, 
apart from Clifford parallelism, these are the only topological parallelisms admitting an automorphism group of dimension 3 or larger.
In particular, we show that a topological parallelism admitting the irreducible 
action of $\SO$ is Clifford. 


\end{abstract}

\section{Introduction}

By a parallelism, we shall always mean an oriented or non-oriented (= ordinary) 
topological parallelism $\Pi$ on $\DP$; 
definitions are given at the end of this introduction.
The automorphism group $\Delta = \mathop{\rm Aut} \Pi$ of a parallelism is a compact Lie subgroup of 
$\mathop{\rm PGL(4,\BR)}$, see \cite {unzush} and \cite{reduzibel}. In particular, $\Delta$ is contained in the maximal compact 
$\mathop{\rm PO(4,\BR)}\le \mathop{\rm PGL(4,\BR)}$.
It is known that only (oriented or non-oriented) Clifford parallelism has $\dim \Delta \ge 4$, see \cite{dim4}, 
\cite{reduzibel} 
(in fact, $\Delta = \mathop{\rm PSO(4,\BR)}$ is 6-dimensional in this case). 

The only possible 3-dimensional group is $\SO$, but it can act in two different ways. 
The reducible action is induced by the natural action of $\SO$ on $\BR \oplus \BR^3$, and the 
irreducible one is induced by $\rm Spin_3 = \mathop{\rm SU}(2,\BC)$ on $\BR^4 = \BC^2$.
Both groups act on Clifford parallelism; the types are represented by the diagonal and by the factors of
$\mathop{\rm PSO(4,\BR)} \cong \SO \times \SO$, respectively. All parallelisms  admitting the reducible $\SO$-action
have been described in \cite {reduzibel}.

Here we show that only Clifford parallelism admits the irreducible $\SO$-action.
This implies the result of \cite{dim4} mentioned above, because every subgroup $\Phi \le
\mathop{\rm PSO}(4,\BR)$ with $\dim \Phi \ge 4$ contains a subgroup $\SO$ with irreducible action.
Non-Clifford examples of regular ordinary parallelisms with a two-dimensional (torus) group
have been given in \cite{dim3}, compare also \cite{gl} for a fresh view of the construction. Regular parallelisms
with an even smaller group are described in \cite{small}. (A topological parallelism is  regular iff its elements are all 
isomorphic to the complex spread.)\\

Let $\CL$ be the line space of $\DP$ and let $\CL^+$ be the twofold covering space of $\CL$ 
consisting of the oriented lines. A topological spread or oriented spread is a compact connected subset 
of $\CL$ or $\CL^+$, respectively, such that each point is incident with exactly one member of $\CC$. 
A (topological) parallelism or oriented parallelism is a set $\Pi$ of topological spreads or oriented spreads, 
respectively, such that every line resp. every oriented line belongs to exacctly one of them, and such that $\Pi$
is compact in the Hausdorff topology on the hyperspace of compact subsets of $\CL$ resp. $\CL^+$. 

For a brief introduction to the basic properties, we refer the reader to \cite{zush}, \cite{unzush} 
and \cite{reduzibel}. A crucial fact is that $\Pi$ is homeomorphic to a star of lines or of oriented lines 
in the ordinary or oriented case, respectively, that is, to the projective plane or to the 2-sphere. 
A homeomorphism is obtained by sending every class to its unique member in the star.
Clifford parallelism consists of the orbits in $\CL$ or $\CL^+$ of one arbitrarily chosen $\SO$-factor of
$\mathop{\rm PSO}(4,\BR)$. The fact that there are two choices will be at the heart of our proof.

\section{Parallelisms with reducible $\SO$-action}

Details of the following construction can be found in \cite{reduzibel}; 
compare also \cite{nonreg} for the non-oriented case. 
The regular (complex) spread $\CC_\BC$ may be described in $\BR^3$ as a the union of a set of reguli carried by one-sheeted 
hyperboloids with common axis of rotation $Z$ and common center $o \in Z$, together with the line at infinity 
defined by $Z^\perp$. The  radius $r$ of the smallest rotation orbit and the asymptotic slope $a$ of the 
hyperboloids have to be related by $a = {c \over r}$, 
where $0<c \in \BR$ is fixed; 
this characterizes the complex spread.
If the ordinary rotation group $\SO$ with fixed pint $o$ is applied to $\CC_\BC$, then one obtains
Clifford parallelism. The same works with orientation, as well. But if the group $\Omega =\SO$ is replaced by a 
conjugate $\Omega^\tau$ where $\tau$ is a nontrivial translation in the direction of $Z$, then the non-oriented 
complex spread does not yield a parallelism; one has to orient $\CC_\BC$,
and in this way one obtains a non-Clifford regular oriented parallelism admitting the group $\Omega^\tau$. 

More generally, the same construction may be applied to a  rotational Betten spread $\CC$, defined by a relaxation of the condition
$a = {c \over r}$, and by allowing the centers of the hyperboloids to be different. If the center is constant and agrees with the 
fixed point of the rotation group, one obtains both ordinary and 
oriented parallelisms. Only the latter occur when the center does not satisfy this condition.  

\section{Results} 

First we recall the left and right Clifford parallelisms on $\DP$, which will play a central role in the proof.
Consider $\BR^4$ as the quaternion skew field $\BH$. Then the orthogonal group $\mathop{\rm SO}(4,\BR)$ may be described as the product of two 
commuting copies $\tilde \Lambda, \tilde \Upsilon$ of the unitary group $\mathop{\mathrm U}(2,\BC)$, consisting 
of the maps $q \mapsto aq$ and $q\mapsto qb$, respectively, 
where $a$, $b$ are 
quaternions of norm one and multiplication is quaternion multiplication. The intersection of the two factors is of order two, 
containing the map $-\rm id$.
Thus, passing to projective space, we get $\mathop{\rm PSO}(4,\BR) = \Lambda \times \Upsilon$, 
a direct product of two copies of $\rm SO (3,\BR)$.
The left and right (ordinary) Clifford parallelisms are defined as the sets of orbits of $\Lambda$ and $\Upsilon$, 
respectively, in the line space $\CL$ of 
$\rm PG(3,\BR)$. For the oriented Clifford parallelisms, one takes the orbits in the space $\CL^+$ of oriented 
lines instead of those in $\CL$. 

The two Clifford parallelisms are equivalent under quaternion conjugation $q \to \bar q$; 
this is immediate from their definition in view of the fact that conjugation does not change the norm and is an anti-automorphism, i.e., 
that $\overline {pq} = \bar q \bar p$. Notice that both $\Lambda$ and $\Upsilon$ are transitive on the point 
set of projective space. Since they centralize
one another, each acts transitively on the parallelism defined by the other, and the group $\mathop{\rm PSO}(4,\BR)$ leaves both parallelisms invariant,
hence it consists of automorphisms of both. 

\bthm\label{main}
Let $\Pi$ be a topological parallelism or oriented parallelism on $\DP$ such that $\dim \mathop{\rm Aut}\Pi \ge 3$.
Then either $\Pi$ is Clifford or $\Pi$ is one of the parallelisms or oriented parallelisms with reducible $\SO$-action 
constructed in {\rm \cite{reduzibel}}, compare Section 2.
\ethm

According to \cite{reduzibel}, Theorem 5.8, the full automorphism of the non-Clifford parallelisms of Section 2 
is the reducible group $\SO$ that was used in order to construct them. Its action on $\Pi$ is transitive.

\bpf
1) It is known that $\dim \mathop{\rm Aut}\Pi \ge 4$ characterizes Clifford parallelism, see \cite{dim4} and,   
for the orientable case, \cite{reduzibel}. However, we do not need this fact. We know that 
$\Delta = \mathop{\rm Aut}\Pi$ is compact; hence, it is a subgroup of the maximal compact 
$\mathop{\rm PO}(4,\BR) \le \mathop{\rm PGL}(4,\BR)$, the identity component of which is a direct product $\Lambda \times \Upsilon$ 
of two copies of $\SO$
(notation as used in the Introduction). $\SO$ is simple and
does not contain any 2-dimensional closed subgroups, hence
a group $\Delta $ of dimension at least 3 either contains one of the factors $\Lambda$ or $\Upsilon$ 
(a rotation group $\SO$ with irreducible action),
or contains a conjugate copy of the diagonal ($\SO$ with reducible action). The latter is a maximal closed connected subgroup, 
hence a group $\Delta$ of dimension $> 3$ 
contains an irreducible $\SO$. All parallelisms with a reducible $\SO$ have been determined in \cite{reduzibel}, hence it suffices to show that 
a parallelism admitting an action of $\Upsilon$ is Clifford. Moreover, it is enough to treat oriented parallelisms, 
because every ordinary parallelism is covered by an oriented one, see \cite{reduzibel}. 
This fact enables a considerable simplification of our arguments.

2) So suppose that $\Upsilon$ acts on an oriented parallelism $\Pi$. We shall make use of the group $\Lambda$ as well, but keep in mind that $\Lambda$ does not act on $\Pi$. Notice first that $\Upsilon$ cannot act trivially on a compact spread, because a group acting trivially on a spread is at most
2-dimensional, compare \cite{CPP}, 61.4.
As both $\Pi$ and every spread $\CC\in \Pi$ 
are homeomorphic to the 2-sphere $\CS_2$, we see that either 
$\Upsilon$ acts trivially  on $\Pi$ and transitively on every member of $\Pi$ (which means that $\Pi$ consists of the $\Upsilon$-orbits in $\CL^+$ and hence coincides with the \it
right \rm Clifford parallelism) or $\Upsilon$ is transitive on $\Pi$. We proceed to show that in the transitive case,
$\Pi$ is the \it left \rm Clifford parallelism. 

3) In \cite{reduzibel}, 2.4, we have shown that the map which sends an oriented line to the pair of its classes with respect to the left and right
oriented Clifford parallelism is a homeomorphism $\CL^+ \to \CS_2\times \CS_2$. The fibers of this product are the left and right Clifford 
parallel classes or, equivalently, the orbits of $\Lambda$ and of $\Upsilon$. We use the homeomorphism to introduce coordinates 
$x,y \in \CS_2$ on $\CL^+$,  writing $(x,y)$ for a line $L \in \CL^+$ whose left and right Clifford classes 
$\Lambda(L)$ and $\Upsilon(L)$ correspond to $x$ and 
$y$, respectively. Then for $\lambda \in \Lambda$ we have
    $$\lambda (x,y) = \lambda(L)= \bigl(\Lambda\lambda(L),\Upsilon\lambda(L)\bigr) 
        = \bigl(\Lambda(L),\lambda\Upsilon(L)\bigr) = (x,\lambda(y)),$$ 
because the groups $\Lambda$ and $\Upsilon$ centralize each other. Likewise, $\upsilon(x,y) = (\upsilon(x),y)$ 
for $\upsilon \in \Upsilon$. 

4) Let now $\CC \in \Pi$ be a parallel class. Then $\Upsilon(\CC) = \Pi$ is a parallelism, hence every element of $\CL^+$ 
belongs to some class $\upsilon(\CC)$, $\upsilon \in \Upsilon$.
This means that $\CC$ meets all orbits of $\Upsilon$. In other words, the projection $p$ from $\CC$ to the second factor of 
$\CS_2\times \CS_2$ is surjective. 

5) We claim that the first projection of $\CC$ is constant, i.e., that $\CC$ is a fiber of the product. In other words, $\CC$ is a $\Lambda$-orbit,
and hence is a class of the left Clifford parallelism. The action of $\Upsilon$ then shows that $\Pi$ is the left Clifford parallelism, 
which completes the proof. 

Let us now prove the claim. 
The stabilizer $\Upsilon_\CC \cong \So$ fixes some pair of elements $x_0, -x_0 \in \CS_2$. If $\CC$ contains a line $L = (x,y)$ such that $x \notin 
\{x_0, -x_0\}$, then $C = \Upsilon_\CC (x)$ is a circle and, by the product action, $C \times \{y\}\subseteq \CC$. There are elements 
$\xi \in \Upsilon \setminus \Upsilon_\CC$ fixing some element $c \in C$.  Then $\CC$ and $\xi(\CC)$ are distinct 
parallel classes, but they have the line $(c,y)$ in common. This is a contradiction. Since $\CC$ is connected, it follows that the first projection 
$\CC \to \{x_0,-x_0\}$ is in fact constant.
\epf

We specialize our result to the case of regular parallelisms.

\bcor\label{cor} 
Let $\Pi$ be a topological parallelism or oriented parallelism on $\DP$ containing a regular spread 
$\CC$ as one of its equivalence classes.
Moreover assume that $\dim \mathop{\rm Aut}\Pi \ge 3$.

a) If $\Pi$ is an ordinary parallelism, then $\Pi$ is (equivalent to) Clifford parallelism.

b) If $\Pi$ is an oriented parallelism, then either $\Pi$ is Clifford or $\Pi$ is one of the examples 
$\Omega^\tau(\CC_\BC)$ given in 
Section 2. In the latter case, $\mathop{\rm Aut}\Pi$ is  $\SO$ with reducible action.
\ecor

\bpf Assertion (a) is obtained by combining Theorem \ref{main} with Corollary 5.7 of \cite{reduzibel}. 
If $\Pi$ is oriented, Theorem \ref{main} tells us 
that $\Pi$ is a parallelism if and only if $\Pi = \Omega^\prime(\CC_\BC)$ 
like in the examples of Section 2, where $\CC_\BC$ is the complex spread and $\Omega'$ is an admissible copy of $\SO$ 
as defined in \cite {reduzibel}, 5.1. According to \cite{reduzibel}, 5.2, this means that 
$\Omega^\prime = \Omega^\mu$ is conjugate
to the rotation group $\Omega = \SO$ with fixed point $o = (0,0,0)$ as in Section 2, via a map of the form 
$\mu(x,y,z) =(x,y,sz+t)$, $s,t \in \BR$, $s \ne 0$. (Here we assume that the line $Z$ in Section 2 is 
given by $x =y=0$). Now since $\CC_\BC$ is the complex spread, the map $(x,y,z) \mapsto (sx,sy,z)$ is an automophism 
of $\CC$, and $(x,y,z) \mapsto (sx,sy,sz)$ centralizes $\Omega$, hence it suffices to use $\mu = \tau$, a 
translation in the direction of $Z$, as in Section 2. According to
\cite{reduzibel}, Theorem 5.8, the automorphism group of $\Pi$ is $\Omega^\tau$.
\epf

%%%%%%%%%%%%%%%%%%%%%%%%%%%%%%%%%%%%%%%%%%%%%%%%%%%%%%%%%%%%%%%%%%%%
%%%%%%%%%%%%%%%%%%%%%%%%%%%%%%%%%%%%%%%%%%%%%%%%%%%%%%%%%%%%%%%%%%%%

\bibliographystyle{plain}

\bigskip
\bigskip
\noindent Rainer L\"owen, Institut f\"ur Analysis und Algebra,
Technische Universit\"at Braunschweig,
Universit\"atsplatz 2,
D 38106 Braunschweig,
Germany

\end{document}